\documentclass[11pt]{article}
\usepackage{amsfonts}
\usepackage{bbm}
\usepackage{mathrsfs}
\usepackage{amscd}
\usepackage{color}
\usepackage{amsmath,amsfonts,amssymb,amscd}
\usepackage{indentfirst,graphics,epsfig,psfrag}
\input{epsf}
\usepackage{ifpdf}
\usepackage{enumerate}
\usepackage{appendix}
\usepackage{enumerate}

\usepackage{lineno}

\makeatletter \@addtoreset{figure}{section} \makeatother
\makeatletter
\long\def\@makecaption#1#2{%
   \vskip 10\p@
   \setbox\@tempboxa\hbox{{#1}\ \ #2}%
   \ifdim \wd\@tempboxa >\hsize

       {#1}\ \ #2\par
   \else
       \hbox to\hsize{\hfil\box\@tempboxa\hfil}%
   \fi}
\makeatother

\begin{document}
\title{\textbf{The resistance distance and Kirchhoff index
on quadrilateral graph and  pentagonal graph}}
\author{
\small  Qun Liu~$^{a,c}$, Zhongzhi Zhang~$^{a,b}$\\
\setcounter{footnote}{-1 }\thanks{Corresponding
author: Qun Liu, E-mail address: liuqun@fudan.edu.cn, zhangzz@fudan.edu.cn}\\
\small  a. School of Computer Science, Fudan University, Shanghai 200433, China\\
\small  b. Shanghai Key Laboratory of Intelligent Information Processing, Fudan University, \\
\small  Shanghai 200433, China\\
\small  c. School of
Mathematics and Statistics, Hexi University,
Zhangye, Gansu, 734000, China.}
\date{}
\maketitle
\begin{abstract}
The quadrilateral graph $Q(G)$ is obtained from $G$
by replacing each edge in $G$ with two parallel paths
of length 1 and 3, whereas the pentagonal graph $W(G)$ is obtained
from $G$ by replacing each edge in $G$ with two parallel paths
of length 1 and 4.
In this paper, closed-form
formulas of resistance distance and Kirchhoff index
for quadrilateral graph and  pentagonal graph are obtained
whenever $G$ is an arbitrary graph.
\\[2mm]

\\[1mm]

\end{abstract}

\section{ Introduction }
All graphs considered in this paper are simple and undirected.
Let $G=(V(G),E(G))$ be a graph with vertex set $V(G)$ and edge set $E(G)$. Let $d_{i}$ be the degree of vertex $i$ in $G$ and $D_{G}=diag(d_{1}, d_{2},\cdots d_{|V(G)|})$ the diagonal matrix with all vertex degrees of $G$ as its diagonal entries. For a graph $G$, let $A_{G}$ and $B_{G}$ denote the adjacency matrix and vertex-edge incidence matrix of $G$, respectively. The matrix $L_{G}=D_{G}-A_{G}$ is called the Laplacian matrix of $G$, where $D_{G}$ is the diagonal matrix of vertex degrees of $G$.

The resistance distance between vertices $u$ and $v$ of $G$ was defined by Klein and Randi$\acute{c}$ $\cite{KR}$ to be the effective resistance between nodes $u$ and $v$ as computed with Ohm's law when all the edges of $G$
are considered to be unit resistors.
The Kirchhoff index $Kf(G)$ was defined in $\cite{KR}$
as $Kf(G)=\sum_{u<v}r_{uv}$, where $r_{uv}(G)$ denotes the resistance distance between $u$ and $v$ in $G$. Resistance distance are, in
fact, intrinsic to the graph, with some nice purely mathematical
interpretations and other interpretations. The Kirchhoff
index was introduced in chemistry as a better alternative
to other parameters used for discriminating different molecules
with similar shapes and structures. See $\cite{KR}$.
The resistance distance
and the Kirchhoff index have attracted extensive attention due to its wide applications in physics, chemistry and others.
Up till now, many results on the resistance distance
and the Kirchhoff index are obtained. See
$(\cite{BG}, \cite{ZH}, \cite{LZhB}-\cite{YK})$ and the references therein to know more.

Recently the quadrilateral graph and the
pentagonal graph have attracted interesting
study. In \cite{LH}, the normalized Laplacian
spectrum of the quadrilateral graph are obtained
and calculate the multiplicative degree-Kirchhoff
index and the spanning trees of the the quadrilateral
graph. Let $G$ be a simple and connected graph with $n$ vertices and $m$ edges. Replacing each edge of $G$ with two parallel paths of lengths 1 and 3 and we call the resulting graph as
quadrilateral graph. We denote by $N$ the total number of vertices and $E$ the total number of edges of $G$. It is clear that $E=4m$, $N=n+2m$. Based on the quadrilateral
graph, we define the new graph. The pentagonal graph of the graph $G$ is
obtained that replacing each edge of $G$ with two parallel paths of lengths 1 and 4. It is clear that $E=5m$, $N=n+3m$.

The rest of this paper is organized as follows. In Section 2, we list some lemmas and give some preliminary results which are used to
prove our main results. In Section 3, we give
the resistance distance and Kirchhoff index of quadrilateral graph
$Q(G)$ whenever $G$ is an arbitrary graph. In Section 4, we give
the resistance distance and Kirchhoff index of pentagonal graph
$W(G)$ whenever $G$ is an arbitrary graph.

\section{Preliminaries}
The $\{1\}$-inverse of $M$ is a matrix $X$ such that $MXM=M$. If $M$ is singular, then it has infinite
$\{1\}$- inverse $\cite{A. Ben-Israel}$.
For a square matrix $M$, the group inverse of $M$, denoted by $M^{\#}$, is the unique matrix $X$ such that $MXM=M$, $XMX=X$ and $MX=XM$. It is known that $M^{\#}$ exists if and only if $rank(M)=rank(M^{2})$ $(\cite{A. Ben-Israel},\cite{BSZHW})$. If $M$ is real symmetric, then $M^{\#}$ exists and $M^{\#}$ is a symmetric $\{1\}$- inverse of $M$. Actually, $M^{\#}$ is equal to the Moore-Penrose inverse of $M$ since $M$ is symmetric $\cite{BSZHW}$.

It is known that resistance distances in a connected graph $G$ can be obtained from any $\{1\}$- inverse of $G$ $(\cite{BG})$. We use $M^{(1)}$ to denote any $\{{1}\}$-inverse of a matrix $M$, and let $(M)_{uv}$ denote the $(u,v)$-entry of $M$.

{\bf Lemma 2.1} $(\cite{BSZHW})$\ Let $G$ be a connected graph. Then
$$
r_{uv}(G)=(L^{(1)}_{G})_{uu}+(L^{(1)}_{G})_{vv}-(L^{(1)}_{G})_{uv}
-(L^{(1)}_{G})_{vu}=(L^{\#}_{G})_{uu}+(L^{\#}_{G})_{vv}-2(L^{\#}_{G})_{uv}.$$

Let $\bf{1_{n}}$ denotes the column vector of dimension $n$ with all the entries equal one. We will often use $\bf{1}$ to denote an all-ones column vector if the dimension can be read from the context.

{\bf Lemma 2.2} $(\cite{BYZhZH})$ \ For any graph $G$, we have
$L^{\#}_{G}\bf{1}$$=0.$

{\bf Lemma 2.3} $(\cite{FZZh})$ \ Let
 \[
\begin{array}{crl}
M=\Bigg(
  \begin{array}{cccccccccccccccc}
   A& B  \\
   C & D \\
 \end{array}
  \Bigg)
\end{array}
\]
be a nonsingular matrix. If $A$ and $D$ are nonsingular,
then
\[
\begin{array}{crl}
M^{-1}&=&\Bigg(
  \begin{array}{cccccccccccccccc}
   A^{-1}+A^{-1}BS^{-1}CA^{-1}&-A^{-1}BS^{-1} \\
   -S^{-1}CA^{-1} & S^{-1}\\
 \end{array}
  \Bigg)
\\&=&\Bigg(
  \begin{array}{cccccccccccccccc}
   (A-BD^{-1}C)^{-1} &-A^{-1}BS^{-1} \\
   -S^{-1}CA^{-1} & S^{-1}\\
 \end{array}
  \Bigg),
\end{array}
\]
where $S=D-CA^{-1}B.$

For a square matrix $M$, let $tr(M)$ denote the trace of $M$.

{\bf Lemma 2.4} $(\cite{SWZhB})$ \ Let $G$ be a connected graph on $n$ vertices. Then
$$Kf(G)=ntr(L^{(1)}_{G})-1^{T}L^{(1)}_{G}1=ntr(L^{\#}_{G}).$$

{\bf Lemma 2.7} \ Let
 \[
\begin{array}{crl}
L=\left(
  \begin{array}{cccccccccccccccc}
   A& B \\
   B^{T} & D \\
 \end{array}
  \right)
\end{array}
\]
be a symmetric block matrix.
If $D$ is nonsingular, then
 \[
\begin{array}{crl}
X=\left(
  \begin{array}{cccccccccccccccc}
   H^{\#}& -H^{\#}BD^{-1} \\
   -D^{-1}B^{T}H^{\#} & D^{-1}+D^{-1}B^{T}H^{\#}BD^{-1} \\
 \end{array}
  \right)
\end{array}
\]
is a symmetric $\{1\}$-inverse of $L$,
where $H=A-BD^{-1}B^{T}$.

\section{The resistance distance and Kirchhoff index of quadrilateral graph}
In this section, we focus on determing the resistance distance and Kirchhoff index of $Q(G)$ whenever $G$ is an arbitrary graph. Let
$E_{G}=\{e_{1},e_{2},\ldots,e_{m}\}$. For each edge $e_{i}
=u_{i}v_{i}\in E_{G}$, there exist two parallel paths of lengths
1 and 3 in $Q(G)$ corresponding to it, which are denoted by
$u_{i}v_{i}$ and $u_{i}u_{i1}u_{i2}v_{i}$ for $i=1,2,\ldots,m$.
Let $V_{1}=\{u_{11},u_{21},\ldots,u_{m1}\}$,
$V_{2}=\{u_{12},u_{22},\ldots,u_{m2}\}$. Then
$V_{Q(G)}=V\bigcup V_{1}\bigcup V_{2}$, where $V$ is the
set of all the vertices inherited from $G$. Our main
results in the following gives the explicit formula of
the resistance distance and Kirchhoff index of $Q(G)$.

{\bf Theorem 3.1} \ Let $G$ be a connected graph with $n$ vertices and $m$ edges and $Q(G)$ be its quadrilateral graph. Then we have the resistance distance
and Kirchhoff index as follows:

(i)For any $i,j\in V(G)$, we have
\begin{eqnarray*}
r_{ij}(Q(G))&=&\frac{3}{4}(L^{\#}_{G})_{ii}+\frac{3}{4}(L^{\#}_{G})_{jj}
-\frac{3}{2}(L^{\#}_{G})_{ij}=\frac{3}{4}r_{ij}(G).
\end{eqnarray*}

(ii)For any $i\in V$, $j\in V_{1}$ or $V_{2}$, $j\in V$, we have
\begin{eqnarray*}
r_{ij}(Q(G))&=&\frac{3}{4}(L^{\#}_{G})_{ii}
+\left[(\frac{1}{2}B^{T}_{1}+
\frac{1}{4}B^{T}_{2})L^{\#}_{G}
(\frac{2}{3}B_{1}+\frac{1}{3}B_{2})\right]_{jj}\\&&
-2\left[L^{\#}_{G}(\frac{1}{2}B_{1}+\frac{1}{4}B_{2})\right]_{ij}.
\end{eqnarray*}

(iii)For any $i\in V_{1}$, $j\in V_{2}$, we have
\begin{eqnarray*}
r_{ij}(Q(G))&=&\frac{4}{3}+\left[(\frac{1}{2}B^{T}_{1}+
\frac{1}{4}B^{T}_{2})L^{\#}_{G}(\frac{2}{3}B_{1}
+\frac{1}{3}B_{2})\right]_{ii}+\left[(\frac{1}{4}B^{T}_{1}+
\frac{1}{2}B^{T}_{2})L^{\#}_{G}
(\frac{1}{3}B_{1}+\frac{2}{3}B_{2})\right]_{jj}\\&&
-[\frac{1}{3}I_{m}+(\frac{1}{2}B^{T}_{1}+
\frac{1}{4}B^{T}_{2})
L^{\#}_{G}(\frac{1}{3}B_{1}+\frac{2}{3}B_{2})]_{ij}.
\end{eqnarray*}

(iv)For any $i,j\in V_{1}$ or $V_{2}$, we have
\begin{eqnarray*}
r_{ij}(Q(G))&=&\frac{4}{3}+
[(\frac{1}{2}B^{T}_{1}+
\frac{1}{4}B^{T}_{2})L^{\#}_{G}(\frac{2}{3}B_{1}
+\frac{1}{3}B_{2})]_{ii}+[(\frac{1}{2}B^{T}_{1}+
\frac{1}{4}B^{T}_{2})L^{\#}_{G}(\frac{2}{3}B_{1}
+\frac{1}{3}B_{2})]_{jj}\\&&-
2[(\frac{1}{2}B^{T}_{1}+
\frac{1}{4}B^{T}_{2})L^{\#}_{G}(\frac{2}{3}B_{1}
+\frac{1}{3}B_{2})]_{ij}.
\end{eqnarray*}

(v)
\begin{eqnarray*}
Kf(Q(G))&=&(n+2m)\left(\frac{3}{4n}Kf(G)
+\frac{5}{12}\left(tr(B^{T}_{1}L^{\#}_{G}B_{1})
+tr(B^{T}_{1}L^{\#}_{G}B_{2})\right)\right.\\&&\left.+\frac{1}{3}
\left(tr(B^{T}_{2}L^{\#}_{G}B_{1})
+tr(B^{T}_{2}L^{\#}_{G}B_{2})\right)\right)
-\frac{3}{4}\left(1^{T}B^{T}_{1}L^{\#}_{G}B_{1}1
+1^{T}B^{T}_{2}L^{\#}_{G}B_{2}1\right.\\&&\left.
+1^{T}B^{T}_{1}L^{\#}_{G}B_{2}1+
1^{T}B^{T}_{2}L^{\#}_{G}B_{1}1\right)-2m,
\end{eqnarray*}
where
$B_{1}+B_{2}=B, B_{1}B^{T}_{2}+B^{T}_{2}B_{1}=A_{G}$ and
$B_{1}B_{1}^{T}+B_{2}B^{T}_{2}=D_{G}$.

{\bf Proof} \ Let $A_{G}$, $D_{G}$ and $B_{G}$ be the adjacency matrix, degree matrix and incidence matrix of $G$. With a suitable labeling for vertices of $Q(G)$,
the Laplacian matrix of $Q(G)$ can be written as follows:
\[
\begin{array}{crl}
 L_{Q_{G}}=\left(
  \begin{array}{cccccccccccccccc}
   2D_{G}-A_{G}& -B_{1}& -B_{2}\\
  -B^{T}_{1} & 2I_{m}& -I_{m} \\
  -B^{T}_{2}& -I_{m}& 2I_{m}\\
 \end{array}
  \right),
\end{array}
\]
where
$B_{1}+B_{2}=B, B_{1}B^{T}_{2}+B^{T}_{2}B_{1}=A_{G}$ and
$B_{1}B_{1}^{T}+B_{2}B^{T}_{2}=D_{G}$.

Let
$A=2D_{G}-A_{G}$,
 $B=\left(
  \begin{array}{ccccc}
 -B_{1}& -B_{2} \\
  \end{array}
\right)$,
 $B^{T}=\left(
  \begin{array}{ccccc}
   -B^{T}_{1}  \\
  -B^{T}_{2}\\
  \end{array}
\right)$ and
$
 D=\left(
  \begin{array}{cccccccccccccccc}
  2I_{m}& -I_{m} \\
   -I_{m}& 2I_{m}\\
 \end{array}
  \right).
$

By Lemma 2.3, it is easily obtained that
$D^{-1}=
 \left(
  \begin{array}{cccccccccccccccc}
    \frac{2}{3}I_{m}&\frac{1}{3}I_{m}\\
  \frac{1}{3}I_{m}&\frac{2}{3}I_{m}\\
 \end{array}
  \right).$

First we begin with the computation of $\{1\}$-inverse of $Q(G)$.

By Lemma 2.7, we have
\[
\begin{array}{crl}
 H&=&2D_{G}-A_{G}-\left(
  \begin{array}{cccccccccccccccc}
  -B_{1}&-B_{2} \\
 \end{array}
  \right)
 \left(
  \begin{array}{cccccccccccccccc}
   \frac{2}{3}I_{m}&\frac{1}{3}I_{m}\\
  \frac{1}{3}I_{m}&\frac{2}{3}I_{m}\\
 \end{array}
  \right)
  \left(
  \begin{array}{cccccccccccccccc}
    -B_{1}^{T}\\
  -B_{2}^{T}\\
 \end{array}
  \right)\\
  &=&2D_{G}-A_{G}-\left(
  \begin{array}{cccccccccccccccc}
    -\frac{2}{3}B_{1}-\frac{1}{3}B_{2}&-\frac{1}{3}B_{1}-\frac{2}{3}B_{2}\\
 \end{array}
  \right)
  \left(
  \begin{array}{cccccccccccccccc}
     -B_{1}^{T} \\
  -B_{2}^{T}\\
   \end{array}
   \right)\\
  &=&2D_{G_{}}-A_{G}-\frac{2}{3}(B_{1}B^{T}_{1}+B_{2}B^{T}_{2})
  -\frac{1}{3}(B_{1}B^{T}_{2}+B_{2}B^{T}_{1})\\
  &=&\frac{4}{3}L^{\#}_{G}
\end{array}
\]
so $H^{\#}=\frac{3}{4}L^{\#}_{G}$.

According to Lemma 2.7, we calculate $-H^{\#}BD^{-1}$ and $-D^{-1}B^{T}H^{\#}$.
\[
\begin{array}{crl}
-H^{\#}BD^{-1}&=&-\frac{3}{4}L^{\#}_{G}\left(
  \begin{array}{cccccccccccccccc}
 -B_{1}& -B_{2}\\
 \end{array}
  \right)
  \left(
  \begin{array}{cccccccccccccccc}
   \frac{2}{3}I_{m}& \frac{1}{3}I_{m}\\
   \frac{1}{3}I_{m}&\frac{2}{3}I_{m}\\
 \end{array}
  \right)\\
  &=&-\frac{3}{4}L^{\#}_{G}\left(
  \begin{array}{cccccccccccccccc}
 -\frac{2}{3}B_{1}-\frac{1}{3}B_{2}&-\frac{1}{3}B_{1}-\frac{2}{3}B_{2}\\
 \end{array}
  \right)\\
   &=&\left(
  \begin{array}{cccccccccccccccc}
 L^{\#}_{G}(\frac{1}{2}B_{1}+\frac{1}{4}B_{2})&
 L^{\#}_{G}(\frac{1}{4}B_{1}+\frac{1}{2}B_{2})\\
 \end{array}
  \right)\\
\end{array}
\]
and
\[
\begin{array}{crl}
-D^{-1}B^{T}H^{\#}=-(H^{\#}BD^{-1})^{T}
  =\left(
  \begin{array}{cccccccccccccccc}
   (\frac{1}{2}B^{T}_{1}+\frac{1}{4}B_{2}^{T})L^{\#}_{G}\\
   (\frac{1}{4}B^{T}_{1}+\frac{1}{2}B_{2}^{T})L^{\#}_{G} \\
 \end{array}
  \right).\\
\end{array}
\]

We are ready to compute $D^{-1}B^{T}H^{\#}BD^{-1}$.
\[
\begin{array}{crl}
D^{-1}B^{T}H^{\#}BD^{-1}
&=&\begin{array}{crl}
 \left(
  \begin{array}{cccccccccccccccc}
     \frac{2}{3}I_{m}& \frac{1}{3}I_{m}\\
   \frac{1}{3}I_{m}&\frac{2}{3}I_{m}\\
 \end{array}
 \right) \left(
  \begin{array}{cccccccccccccccc}
-B^{T}_{1}  \\
  -B^{T}_{2}\\
 \end{array}
 \right)\frac{3}{4}L^{\#}_{G}
 \left(
  \begin{array}{cccccccccccccccc}
   -B_{1}&-B_{2} \\
 \end{array}
 \right)
 \left(
  \begin{array}{cccccccccccccccc}
     \frac{2}{3}I_{m}& \frac{1}{3}I_{m}\\
   \frac{1}{3}I_{m}&\frac{2}{3}I_{m}\\
 \end{array}
 \right)
\end{array}\\
&=&\left(
\begin{array}{cccccccccccccccc}
(\frac{1}{2}B^{T}_{1}+\frac{1}{4}B^{T}_{2})L^{\#}_{G}(\frac{2}{3}B_{1}+\frac{1}{3}B_{2})&
(\frac{1}{2}B^{T}_{1}+\frac{1}{4}B^{T}_{2})L^{\#}_{G}(\frac{1}{3}B_{1}+\frac{2}{3}B_{2})\\
 (\frac{1}{4}B^{T}_{1}+\frac{1}{2}B^{T}_{2})L^{\#}_{G}(\frac{2}{3}B_{1}+\frac{1}{3}B_{2})&
(\frac{1}{4}B^{T}_{1}+\frac{1}{2}B^{T}_{2})L^{\#}_{G}(\frac{1}{3}B_{1}+\frac{2}{3}B_{2})\\
 \end{array}
 \right).
\end{array}
\]

Let $P=(\frac{1}{2}B^{T}_{1}+
\frac{1}{4}B^{T}_{2})L^{\#}_{G}(\frac{2}{3}B_{1}+\frac{1}{3}B_{2})$,
$Q=(\frac{1}{2}B^{T}_{1}+
\frac{1}{4}B^{T}_{2})L^{\#}_{G}(\frac{1}{3}B_{1}+\frac{2}{3}B_{2})$,
$M=(\frac{1}{4}B^{T}_{1}+
\frac{1}{2}B^{T}_{2})L^{\#}_{G}(\frac{2}{3}B_{1}+\frac{1}{3}B_{2})$
and $N=(\frac{1}{4}B^{T}_{1}+
\frac{1}{2}B^{T}_{2})L^{\#}_{G}(\frac{1}{3}B_{1}+\frac{2}{3}B_{2})$.
Based on Lemma 2.3 and 2.7, the following matrix
\begin{eqnarray}
N=\left(
  \begin{array}{cccccccccccccccc}
\frac{3}{4}L^{\#}_{G}&
L^{\#}_{G}(\frac{1}{2}B_{1}+\frac{1}{4}B_{2})&
 L^{\#}_{G}(\frac{1}{4}B_{1}+\frac{1}{2}B_{2})
 \\
L^{\#}_{G}(\frac{1}{2}B^{T}_{1}+\frac{1}{4}B^{T}_{2})&
\frac{2}{3}I_{m}+P&
\frac{1}{3}I_{m}+Q\\
L^{\#}_{G}(\frac{1}{4}B^{T}_{1}+\frac{1}{2}B^{T}_{2})&
\frac{1}{3}I_{m}+M
&\frac{2}{3}I_{m}+N\\
 \end{array}
  \right)
\end{eqnarray}
is a symmetric $\{1\}$-inverse of $L_{Q(G)}$.

For any $i,j\in V(G)$, by Lemma 2.1 and the Equation $(3.1)$, we have
\begin{eqnarray*}
r_{ij}(Q(G))&=&\frac{3}{4}(L^{\#}_{G})_{ii}+\frac{3}{4}(L^{\#}_{G})_{jj}
-\frac{3}{2}(L^{\#}_{G})_{ij}=\frac{3}{4}r_{ij}(G),
\end{eqnarray*}
as stated in $(i)$.

For any $i\in V$, $j\in V_{1}$ or $V_{2}$, $j\in V$, by Lemma 2.1 and the Equation $(3.1)$, we have
\begin{eqnarray*}
r_{ij}(Q(G))&=&\frac{3}{4}(L^{\#}_{G})_{ii}
+\left[\frac{2}{3}I_{m}+(\frac{1}{2}B^{T}_{1}+
\frac{1}{4}B^{T}_{2})L^{\#}_{G}
(\frac{2}{3}B_{1}+\frac{1}{3}B_{2})\right]_{jj}\\&&
-2\left[L^{\#}_{G}(\frac{1}{2}B_{1}+\frac{1}{4}B_{2})\right]_{ij},
\end{eqnarray*}
as stated in $(ii)$.

For any $i\in V_{1}$, $j\in V_{2}$, by Lemma 2.1 and the Equation $(3.1)$, we have
\begin{eqnarray*}
r_{ij}(Q(G))&=&\frac{4}{3}+\left[(\frac{1}{2}B^{T}_{1}+
\frac{1}{4}B^{T}_{2})L^{\#}_{G}(\frac{2}{3}B_{1}
+\frac{1}{3}B_{2})\right]_{ii}+\left[(\frac{1}{4}B^{T}_{1}+
\frac{1}{2}B^{T}_{2})L^{\#}_{G}
(\frac{1}{3}B_{1}+\frac{2}{3}B_{2})\right]_{jj}\\&&
-[\frac{1}{3}I_{m}+(\frac{1}{2}B^{T}_{1}+
\frac{1}{4}B^{T}_{2})
L^{\#}_{G}(\frac{1}{3}B_{1}+\frac{2}{3}B_{2})]_{ij},
\end{eqnarray*}
as stated in $(iii)$.

For any $i,j\in V_{1}$ or $V_{2}$, by Lemma 2.1 and the Equation $(3.1)$, we have
\begin{eqnarray*}
r_{ij}(Q(G))&=&\frac{4}{3}+
[(\frac{1}{2}B^{T}_{1}+
\frac{1}{4}B^{T}_{2})L^{\#}_{G}(\frac{2}{3}B_{1}
+\frac{1}{3}B_{2})]_{ii}+[(\frac{1}{2}B^{T}_{1}+
\frac{1}{4}B^{T}_{2})L^{\#}_{G}(\frac{2}{3}B_{1}
+\frac{1}{3}B_{2})]_{jj}\\&&-
2[(\frac{1}{2}B^{T}_{1}+
\frac{1}{4}B^{T}_{2})L^{\#}_{G}(\frac{2}{3}B_{1}
+\frac{1}{3}B_{2})]_{ij},
\end{eqnarray*}
as stated in $(iv)$.

By Lemma 2.4, we have
\begin{eqnarray*}
Kf(Q(G))
&=&(n+2m)tr( N)-{\bf{1}^{T}}N{\bf{1}}
\end{eqnarray*}
\begin{eqnarray*}
&=&(n+2m)\left(\frac{3}{4}tr(L^{\#}_{G})
+tr\left[(\frac{1}{2}B_{1}^{T}+\frac{1}{4}B_{2}^{T})
L^{\#}_{G}(\frac{2}{3}B_{1}+\frac{1}{3}B_{2})\right]\right.\\
&&\left.+
tr\left[(\frac{1}{4}B_{1}^{T}+\frac{1}{2}B_{2}^{T})
L^{\#}_{G}(\frac{1}{3}B_{1}+\frac{2}{3}B_{2})\right]
+2tr(\frac{2}{3}I_{m})\right)-
{\bf{1}^{T}}N{\bf{1}}\\
&=&(n+2m)\left(\frac{3}{4n}Kf(G)
+tr\left[(\frac{1}{2}B_{1}^{T}+\frac{1}{4}B_{2}^{T})
L^{\#}_{G}(\frac{2}{3}B_{1}+\frac{1}{3}B_{2})\right]\right.\\
&&\left.+
tr\left[(\frac{1}{4}B_{1}^{T}+\frac{1}{2}B_{2}^{T})
L^{\#}_{G}(\frac{1}{3}B_{1}+\frac{2}{3}B_{2})\right]
+\frac{4m}{3}\right)-
{\bf{1}^{T}}N{\bf{1}}\\
&=&(n+2m)\left(\frac{3}{4n}Kf(G)
+\frac{5}{12}\left(tr(B^{T}_{1}L^{\#}_{G}B_{1})
+tr(B^{T}_{1}L^{\#}_{G}B_{2})\right)\right.\\&&\left.+\frac{1}{3}
\left(tr(B^{T}_{2}L^{\#}_{G}B_{1})
+tr(B^{T}_{2}L^{\#}_{G}B_{2})\right)\right)-
{\bf{1}^{T}}N{\bf{1}}.
\end{eqnarray*}

Since $L^{\#}_{G}1=0$, then

${\bf{1}^{T}}N{\bf{1}^{T}}$
\begin{eqnarray*}
&=&1^{T}[(\frac{1}{2}B^{t}_{1}+
\frac{1}{4}B^{T}_{2})L^{\#}_{G}(\frac{2}{3}B_{1}+\frac{1}{3}B_{2})]1
+1^{T}(\frac{1}{2}B^{t}_{1}+
\frac{1}{4}B^{T}_{2})L^{\#}_{G}(\frac{1}{3}B_{1}+\frac{2}{3}B_{2})1
+2\cdot\frac{1}{3}1^{T}1\\
&&+1^{T}(\frac{1}{4}B^{t}_{1}+
\frac{1}{2}B^{T}_{2})L^{\#}_{G}(\frac{2}{3}B_{1}+\frac{1}{3}B_{2})1
+1^{T}(\frac{1}{4}B^{t}_{1}+
\frac{1}{2}B^{T}_{2})L^{\#}_{G}(\frac{1}{3}B_{1}+\frac{2}{3}B_{2})1
+2\cdot\frac{2}{3}1^{T}1\\
&=&
\frac{3}{4}1^{T}B^{T}_{1}L^{\#}_{G}B_{1}1+\frac{3}{4}1^{T}B^{T}_{2}L^{\#}_{G}B_{2}1
+\frac{3}{4}1^{T}B^{T}_{1}L^{\#}_{G}B_{2}1+
\frac{3}{4}1^{T}B^{T}_{2}L^{\#}_{G}B_{1}1+2m.
\end{eqnarray*}

Plugging the above equation into $Kf(Q(G))$,
we obtain the required result in $(v)$.

\section{The resistance distance and Kirchhoff index of pentagonal graph}
In this section, we focus on determing the resistance distance and Kirchhoff index of $W(G)$ whenever $G$ be an arbitrary graph.
Let $E_{G}=\{e_{1},e_{2},\ldots,e_{m}\}$. For each edge $e_{i}
=u_{i}v_{i}\in E_{G}$, there exist two parallel paths of lengths
1 and 4 in $W(G)$ corresponding to it, which are denoted by
$u_{i}v_{i}$ and $u_{i}u_{i1}u_{i2}u_{i3}v_{i}$ for $i=1,2,\ldots,m$.
Let $V_{1}=\{u_{11},u_{21},\ldots,u_{m1}\}$,
$V_{2}=\{u_{12},u_{22},\ldots,u_{m2}\}$,
$V_{3}=\{u_{13},u_{23},\ldots,u_{m3}\}$. Then
$V_{W(G)}=V\bigcup V_{1}\bigcup V_{2}\bigcup V_{3}$, where $V$ is the
set of all the vertices inherited from $G$. Our main
results in the following gives the explicit formula of
the resistance distance and Kirchhoff index of $W(G)$.

{\bf Theorem 4.1} \ Let $G$ be a connected graph with $n$ vertices and $m$ edges and $W(G)$ be its pentagonal graph of $G$.
Then we have the resistance distance
and Kirchhoff index as follows:

(i)For any $i,j\in V(G)$, we have
\begin{eqnarray*}
r_{ij}(W(G))&=&\frac{4}{5}(L^{\#}_{G})_{ii}+\frac{4}{5}(L^{\#}_{G})_{jj}
-\frac{8}{5}(L^{\#}_{G})_{ij}=\frac{4}{5}r_{ij}(G).
\end{eqnarray*}

(ii)For any $i\in V$, $j\in V_{1}$, we have
\begin{eqnarray*}
r_{ij}(W(G))&=&\frac{4}{5}(L^{\#}_{G})_{ii}
+\left[\frac{3}{4}I_{m}+(\frac{3}{5}B^{T}_{1}
+\frac{1}{5}B^{T}_{2})L^{\#}_{G}(\frac{3}{4}B_{1}
+\frac{1}{4}B_{2})\right]_{jj}\\&&
-2\left[L^{\#}_{G}(\frac{3}{5}B_{1}+\frac{1}{5}B_{2})\right]_{ij}.
\end{eqnarray*}

(iii)For any $i\in V$, $j\in V_{2}$, we have
\begin{eqnarray*}
r_{ij}(W(G))&=&\frac{4}{5}(L^{\#}_{G})_{ii}
+\left[I_{m}+(\frac{3}{5}B^{T}_{1}
+\frac{1}{5}B^{T}_{2})L^{\#}_{G}(\frac{3}{4}B_{1}
+\frac{1}{4}B_{2})\right]_{jj}\\&&
-2\left[L^{\#}_{G}(\frac{2}{5}B_{1}+\frac{2}{5}B_{2})\right]_{ij}.
\end{eqnarray*}

(iv)For any $i\in V$, $j\in V_{3}$, we have
\begin{eqnarray*}
r_{ij}(W(G))&=&\frac{4}{5}(L^{\#}_{G})_{ii}
+\left[(\frac{1}{5}B^{T}_{1}+\frac{3}{5}B^{T}_{2})L^{\#}_{G}
(\frac{1}{4}B_{1}+\frac{3}{4}B_{2})\right]_{jj}\\&&
-2\left[L^{\#}_{G}(\frac{1}{5}B_{1}+\frac{3}{5}B_{2})\right]_{ij}.
\end{eqnarray*}

(v)For any $i\in V_{1}$, $j\in V_{2}$, we have
\begin{eqnarray*}
r_{ij}(W(G))&=&\frac{5}{4}+\left[(\frac{3}{5}B^{T}_{1}
+\frac{1}{5}B^{T}_{2})L^{\#}_{G}(\frac{3}{4}B_{1}
+\frac{1}{4}B_{2})\right]_{ii}+\left[(\frac{2}{5}B^{T}_{1}+\frac{2}{5}B^{T}_{2})L^{\#}_{G}
(\frac{1}{2}B_{1}+\frac{1}{2}B_{2})\right]_{jj}\\&&
-[\frac{1}{2}I_{m}+(\frac{3}{5}B^{T}_{1}+\frac{1}{5}B^{T}_{2})
L^{\#}_{G}(\frac{1}{2}B_{1}+\frac{1}{2}B_{2})]_{ij}.
\end{eqnarray*}

(vi)For any $i\in V_{1}$, $j\in V_{3}$, we have
\begin{eqnarray*}
r_{ij}(W(G))&=&\frac{3}{2}+\left[(\frac{3}{5}B^{T}_{1}
+\frac{1}{5}B^{T}_{2})L^{\#}_{G}(\frac{3}{4}B_{1}
+\frac{1}{4}B_{2})\right]_{ii}+\left[(\frac{1}{5}B^{T}_{1}+\frac{3}{5}B^{T}_{2})L^{\#}_{G}
(\frac{1}{4}B_{1}+\frac{3}{4}B_{2})\right]_{jj}\\&&
-[\frac{1}{4}I_{m}+(\frac{3}{5}B^{T}_{1}+\frac{1}{5}B^{T}_{2})
L^{\#}_{G}(\frac{1}{4}B_{1}+\frac{3}{4}B_{2})]_{ij}.
\end{eqnarray*}

(vii)For any $i\in V_{2}$, $j\in V_{3}$, we have
\begin{eqnarray*}
r_{ij}(W(G))&=&\frac{7}{4}+\left[(\frac{2}{5}B^{T}_{1}+\frac{2}{5}B^{T}_{2})L^{\#}_{G}
(\frac{1}{2}B_{1}+\frac{1}{2}B_{2})\right]_{ii}+\left[(\frac{1}{5}B^{T}_{1}+\frac{3}{5}B^{T}_{2})L^{\#}_{G}
(\frac{1}{4}B_{1}+\frac{3}{4}B_{2})\right]_{jj}\\&&
-[\frac{1}{2}I_{m}+(\frac{1}{4}B^{T}_{1}+\frac{1}{4}B^{T}_{2})L^{\#}_{G}
(\frac{1}{4}B_{1}+\frac{3}{4}B_{2})]_{ij}.
\end{eqnarray*}

(viii)
$Kf(W(G))$
\begin{eqnarray*}
&=&(n+3m)\left(\frac{4}{5n}Kf(G)
+\frac{61}{100}\left(tr(B^{T}_{1}L^{\#}_{G}B_{1})
+tr(B^{T}_{2}L^{\#}_{G}B_{2})\right)\right.\\&&\left.+\frac{1}{2}
\left(tr(B^{T}_{1}L^{\#}_{G}B_{2})
+tr(B^{T}_{1}L^{\#}_{G}B_{2})\right)+\frac{5m}{2}\right)-
\frac{141}{80}1^{T}B^{T}_{1}L^{\#}_{G}B_{1}1-
\frac{131}{80}1^{T}B^{T}_{1}L^{\#}_{G}B_{2}1\\&&
-\frac{133}{80}1^{T}B^{T}_{2}L^{\#}_{G}B_{1}1-
\frac{127}{80}1^{T}B^{T}_{2}L^{\#}_{G}B_{2}1-5m,
\end{eqnarray*}
where
$B_{1}+B_{2}=B, B_{1}B^{T}_{2}+B^{T}_{2}B_{1}=A_{G}$ and
$B_{1}B_{1}^{T}+B_{2}B^{T}_{2}=D_{G}$.

{\bf Proof} \ Let $A_{G}$, $D_{G}$ and $B_{G}$ be the adjacency matrix, degree matrix and incidence matrix of $G$. With a suitable labeling for vertices of $W(G)$,
the Laplacian matrix of $W(G)$ can be written as follows:
\[
\begin{array}{crl}
 L_{W(G)}=\left(
  \begin{array}{cccccccccccccccc}
   2D_{G}-A_{G}& -B_{1}& 0&-B_{2}\\
  -B^{T}_{1} & 2I_{m}& -I_{m}&0 \\
  0&-I_{m}&2I_{m}&-I_{m}\\
  -B^{T}_{2}& 0&-I_{m}& 2I_{m}\\
 \end{array}
  \right),
\end{array}
\]
where
$B_{1}+B_{2}=B, B_{1}B^{T}_{2}+B^{T}_{2}B_{1}=A_{G}$ and
$B_{1}B_{1}^{T}+B_{2}B^{T}_{2}=D_{G}$.

Let
$A=2D_{G}-A_{G}$,
 $B=\left(
  \begin{array}{ccccc}
 -B_{1}& 0&-B_{2} \\
  \end{array}
\right)$,
 $B^{T}=\left(
  \begin{array}{ccccc}
   -B^{T}_{1}  \\
   0\\
  -B^{T}_{2}\\
  \end{array}
\right)$ and
$
 D=\left(
  \begin{array}{cccccccccccccccc}
  2I_{m}& -I_{m}&0 \\
   -I_{m}& 2I_{m}&-I_{m}\\
   0&-I_{m}&2I_{m}
 \end{array}
  \right)
$

By Lemma 2.3, we have
$D^{-1}=
 \left(
  \begin{array}{cccccccccccccccc}
    \frac{3}{4}I_{m}&\frac{1}{2}I_{m}&\frac{1}{4}I_{m}\\
  \frac{1}{2}I_{m}&I_{m}&\frac{1}{2}I_{m}\\
  \frac{1}{4}I_{m}&\frac{1}{2}I_{m}&\frac{3}{4}I_{m}\\
 \end{array}
  \right).$

First we begin with the computation of $\{1\}$-inverse of $W(G)$.

By Lemma 2.7, we have
\[
\begin{array}{crl}
 H&=&2D_{G}-A_{G}-\left(
  \begin{array}{cccccccccccccccc}
  -B_{1}&0&-B_{2} \\
 \end{array}
  \right)
 \left(
  \begin{array}{cccccccccccccccc}
  \frac{3}{4}I_{m}&\frac{1}{2}I_{m}&\frac{1}{4}I_{m}\\
  \frac{1}{2}I_{m}&I_{m}&\frac{1}{2}I_{m}\\
  \frac{1}{4}I_{m}&\frac{1}{2}I_{m}&\frac{3}{4}I_{m}\\
 \end{array}
  \right)
  \left(
  \begin{array}{cccccccccccccccc}
    -B_{1}^{T}\\
    0\\
  -B_{2}^{T}\\
 \end{array}
  \right)\\
  &=&2D_{G}-A_{G}-\left(
  \begin{array}{cccccccccccccccc}
    -\frac{3}{4}B_{1}-\frac{1}{4}B_{2}&
    -\frac{1}{2}B_{1}-\frac{1}{2}B_{2}&
    -\frac{1}{4}B_{1}-\frac{3}{4}B_{2}\\
 \end{array}
  \right)
  \left(
  \begin{array}{cccccccccccccccc}
     -B_{1}^{T} \\
     0\\
  -B_{2}^{T}\\
   \end{array}
   \right)\\
  &=&2D_{G_{}}-A_{G}-\frac{3}{4}(B_{1}B^{T}_{1}+B_{2}B^{T}_{2})
  -\frac{1}{4}(B_{1}B^{T}_{2}+B_{2}B^{T}_{1})\\
  &=&\frac{5}{4}L^{\#}_{G},
\end{array}
\]
so $H^{\#}=\frac{4}{5}L^{\#}_{G}$.

According to Lemma 2.7, we calculate $-H^{\#}BD^{-1}$ and $-D^{-1}B^{T}H^{\#}$.
\[
\begin{array}{crl}
-H^{\#}BD^{-1}&=&-\frac{4}{5}L^{\#}_{G}\left(
  \begin{array}{cccccccccccccccc}
 -B_{1}& 0&-B_{2}\\
 \end{array}
  \right)
  \left(
  \begin{array}{cccccccccccccccc}
   \frac{3}{4}I_{m}&\frac{1}{2}I_{m}&\frac{1}{4}I_{m}\\
  \frac{1}{2}I_{m}&I_{m}&\frac{1}{2}I_{m}\\
  \frac{1}{4}I_{m}&\frac{1}{2}I_{m}&\frac{3}{4}I_{m}\\
 \end{array}
  \right)\\
  &=&\frac{4}{5}L^{\#}_{G}\left(
  \begin{array}{cccccccccccccccc}
 \frac{3}{4}B_{1}+\frac{1}{4}B_{2}&
 \frac{1}{2}B_{1}+\frac{1}{2}B_{2}&
 \frac{1}{4}B_{1}+\frac{3}{4}B_{2}\\
 \end{array}
  \right)\\
   &=&\left(
  \begin{array}{cccccccccccccccc}
 L^{\#}_{G}(\frac{3}{5}B_{1}+\frac{1}{5}B_{2})&
 L^{\#}_{G}(\frac{2}{5}B_{1}+\frac{2}{5}B_{2})&
 L^{\#}_{G}(\frac{1}{5}B_{1}+\frac{3}{5}B_{2})\\
 \end{array}
  \right)\\
\end{array}
\]
and
\[
\begin{array}{crl}
-D^{-1}B^{T}H^{\#}=-(H^{\#}BD^{-1})^{T}
  =\left(
  \begin{array}{cccccccccccccccc}
   (\frac{3}{5}B^{T}_{1}+\frac{1}{5}B^{T}_{2})L^{\#}_{G}\\
   (\frac{2}{5}B^{T}_{1}+\frac{2}{5}B^{T}_{2})L^{\#}_{G} \\
  (\frac{1}{5}B^{T}_{1}+\frac{3}{5}B^{T}_{2}) L^{\#}_{G}\\
 \end{array}
  \right).\\
\end{array}
\]

We are ready to compute the $D^{-1}B^{T}H^{\#}BD^{-1}$.

$D^{-1}B^{T}H^{\#}BD^{-1}$
\[
\begin{array}{crl}
&=&\begin{array}{crl}
 \left(
  \begin{array}{cccccccccccccccc}
    \frac{3}{4}I_{m}&\frac{1}{2}I_{m}&\frac{1}{4}I_{m}\\
  \frac{1}{2}I_{m}&I_{m}&\frac{1}{2}I_{m}\\
  \frac{1}{4}I_{m}&\frac{1}{2}I_{m}&\frac{3}{4}I_{m}\\
 \end{array}
 \right) \left(
  \begin{array}{cccccccccccccccc}
-B^{T}_{1}  \\
   0\\
  -B^{T}_{2}\\
 \end{array}
 \right)\frac{4}{5}L^{\#}_{G}
 \left(
  \begin{array}{cccccccccccccccc}
   -B_{1}&0&-B_{2} \\
 \end{array}
 \right)
 \left(
  \begin{array}{cccccccccccccccc}
     \frac{3}{4}I_{m}&\frac{1}{2}I_{m}&\frac{1}{4}I_{m}\\
  \frac{1}{2}I_{m}&I_{m}&\frac{1}{2}I_{m}\\
  \frac{1}{4}I_{m}&\frac{1}{2}I_{m}&\frac{3}{4}I_{m}\\
 \end{array}
 \right)
\end{array}
\end{array}\\
\]

\[
\begin{array}{crl}
&=&\left(
\begin{array}{cccccccccccccccc}
(\frac{3}{5}B^{T}_{1}+\frac{1}{5}B^{T}_{2})L^{\#}_{G}(\frac{3}{4}B_{1}+\frac{1}{4}B_{2})&
(\frac{3}{5}B^{T}_{1}+\frac{1}{5}B^{T}_{2})L^{\#}_{G}(\frac{1}{2}B_{1}+\frac{1}{2}B_{2})&
(\frac{3}{5}B^{T}_{1}+\frac{1}{5}B^{T}_{2})
L^{\#}_{G}(\frac{1}{4}B_{1}+\frac{3}{4}B_{2})
\\
 (\frac{2}{5}B^{T}_{1}+\frac{2}{5}B^{T}_{2})L^{\#}_{G}(\frac{3}{4}B_{1}+\frac{1}{4}B_{2})&
(\frac{2}{5}B^{T}_{1}+\frac{2}{5}B^{T}_{2})L^{\#}_{G}(\frac{1}{2}B_{1}+\frac{1}{2}B_{2})&
(\frac{1}{5}B^{T}_{1}+\frac{1}{2}B^{T}_{2})L^{\#}_{G}(\frac{1}{4}B_{1}+\frac{3}{4}B_{2})\\
(\frac{1}{5}B^{T}_{1}+\frac{3}{5}B^{T}_{2})L^{\#}_{G}(\frac{3}{4}B_{1}+\frac{1}{4}B_{2})&
(\frac{1}{5}B^{T}_{1}+\frac{3}{5}B^{T}_{2})L^{\#}_{G}(\frac{1}{2}B_{1}+\frac{1}{2}B_{2})&
(\frac{1}{5}B^{T}_{1}+\frac{3}{5}B^{T}_{2})L^{\#}_{G}(\frac{1}{4}B_{1}+\frac{3}{4}B_{2})\\
 \end{array}
 \right).
\end{array}\\
\]

Let $P_{1}=(\frac{3}{5}B^{T}_{1}
+\frac{1}{5}B^{T}_{2})L^{\#}_{G}(\frac{3}{4}B_{1}
+\frac{1}{4}B_{2})$, $P_{2}=(\frac{3}{5}B^{T}_{1}+\frac{1}{5}B^{T}_{2})
L^{\#}_{G}(\frac{1}{2}B_{1}+\frac{1}{2}B_{2})$
$P_{3}=(\frac{3}{5}B^{T}_{1}+\frac{1}{5}B^{T}_{2})
L^{\#}_{G}(\frac{1}{4}B_{1}+\frac{3}{4}B_{2})$,
$Q_{1}=(\frac{2}{5}B^{T}_{1}+\frac{2}{5}B^{T}_{2})L^{\#}_{G}
(\frac{3}{4}B_{1}+\frac{1}{4}B_{2})$,
$Q_{2}=(\frac{2}{5}B^{T}_{1}+\frac{2}{5}B^{T}_{2})L^{\#}_{G}
(\frac{1}{2}B_{1}+\frac{1}{2}B_{2})$,
$Q_{3}=(\frac{1}{4}B^{T}_{1}+\frac{1}{4}B^{T}_{2})L^{\#}_{G}
(\frac{1}{4}B_{1}+\frac{3}{4}B_{2})$,
$R_{1}=(\frac{1}{5}B^{T}_{1}+\frac{3}{5}B^{T}_{2})L^{\#}_{G}
(\frac{3}{4}B_{1}+\frac{1}{4}B_{2})$,
$R_{2}=(\frac{1}{5}B^{T}_{1}+\frac{3}{5}B^{T}_{2})L^{\#}_{G}
(\frac{1}{2}B_{1}+\frac{1}{2}B_{2})$,
$R_{3}=(\frac{1}{5}B^{T}_{1}+\frac{3}{5}B^{T}_{2})L^{\#}_{G}
(\frac{1}{4}B_{1}+\frac{3}{4}B_{2})$,

Based on Lemma 2.3 and 2.7, the following matrix
\begin{eqnarray}
N=\left(
  \begin{array}{cccccccccccccccc}
\frac{4}{5}L^{\#}_{G}&L^{\#}_{G}(\frac{3}{5}B_{1}+\frac{1}{5}B_{2})&
 L^{\#}_{G}(\frac{2}{5}B_{1}+\frac{2}{5}B_{2})&
 L^{\#}_{G}(\frac{1}{5}B_{1}+\frac{3}{5}B_{2})\\
(\frac{3}{5}B^{T}_{1}+\frac{1}{5}B^{T}_{2})L^{\#}_{G}
&\frac{3}{4}I_{m}+P_{1}&\frac{1}{2}I_{m}+P_{2}
&\frac{1}{4}I_{m}+P_{3}\\
(\frac{2}{5}B^{T}_{1}+\frac{2}{5}B^{T}_{2})L^{\#}_{G}&
\frac{1}{2}I_{m}+Q_{1}&I_{m}+Q_{2}&\frac{1}{2}I_{m}+Q_{3}\\
 (\frac{1}{5}B^{T}_{1}+\frac{3}{5}B^{T}_{2}) L^{\#}_{G}&\frac{1}{4}I_{m}+R_{1}&\frac{1}{2}I_{m}+R_{2}&
 \frac{3}{4}I_{m}+R_{3}\\
 \end{array}
  \right)
\end{eqnarray}
is a symmetric $\{1\}$-inverse of $L$.

For any $i,j\in V(G)$, by Lemma 2.1 and the Equation $(4.2)$, we have
\begin{eqnarray*}
r_{ij}(W(G))&=&\frac{4}{5}(L^{\#}_{G})_{ii}+\frac{4}{5}(L^{\#}_{G})_{jj}
-\frac{8}{5}(L^{\#}_{G})_{ij}=\frac{4}{5}r_{ij}(G),
\end{eqnarray*}
as stated in $(i)$.

For any $i\in V$, $j\in V_{1}$, by Lemma 2.1 and the Equation $(4.2)$, we have
\begin{eqnarray*}
r_{ij}(W(G))&=&\frac{4}{5}(L^{\#}_{G})_{ii}
+\left[\frac{3}{4}I_{m}+(\frac{3}{5}B^{T}_{1}
+\frac{1}{5}B^{T}_{2})L^{\#}_{G}(\frac{3}{4}B_{1}
+\frac{1}{4}B_{2})\right]_{jj}\\&&
-2\left[L^{\#}_{G}(\frac{3}{5}B_{1}+\frac{1}{5}B_{2})\right]_{ij},
\end{eqnarray*}
as stated in $(ii)$.

For any $i\in V$, $j\in V_{2}$, by Lemma 2.1 and the Equation $(4.2)$, we have
\begin{eqnarray*}
r_{ij}(W(G))&=&\frac{4}{5}(L^{\#}_{G})_{ii}
+\left[I_{m}+(\frac{3}{5}B^{T}_{1}
+\frac{1}{5}B^{T}_{2})L^{\#}_{G}(\frac{3}{4}B_{1}
+\frac{1}{4}B_{2})\right]_{jj}\\&&
-2\left[L^{\#}_{G}(\frac{2}{5}B_{1}+\frac{2}{5}B_{2})\right]_{ij},
\end{eqnarray*}
as stated in $(iii)$.

For any $i\in V$, $j\in V_{3}$, by Lemma 2.1 and the Equation $(4.2)$, we have
\begin{eqnarray*}
r_{ij}(W(G))&=&\frac{4}{5}(L^{\#}_{G})_{ii}
+\left[(\frac{1}{5}B^{T}_{1}+\frac{3}{5}B^{T}_{2})L^{\#}_{G}
(\frac{1}{4}B_{1}+\frac{3}{4}B_{2})\right]_{jj}\\&&
-2\left[L^{\#}_{G}(\frac{1}{5}B_{1}+\frac{3}{5}B_{2})\right]_{ij},
\end{eqnarray*}
as stated in $(iv)$.

For any $i\in V_{1}$, $j\in V_{2}$, by Lemma 2.1 and the Equation $(3.1)$, we have
\begin{eqnarray*}
r_{ij}(W(G))&=&\frac{5}{4}+\left[(\frac{3}{5}B^{T}_{1}
+\frac{1}{5}B^{T}_{2})L^{\#}_{G}(\frac{3}{4}B_{1}
+\frac{1}{4}B_{2})\right]_{ii}+\left[(\frac{2}{5}B^{T}_{1}+\frac{2}{5}B^{T}_{2})L^{\#}_{G}
(\frac{1}{2}B_{1}+\frac{1}{2}B_{2})\right]_{jj}\\&&
-[\frac{1}{2}I_{m}+(\frac{3}{5}B^{T}_{1}+\frac{1}{5}B^{T}_{2})
L^{\#}_{G}(\frac{1}{2}B_{1}+\frac{1}{2}B_{2})]_{ij},
\end{eqnarray*}
as stated in $(v)$.

For any $i\in V_{1}$, $j\in V_{3}$, by Lemma 2.1 and the Equation $(4.2)$, we have
\begin{eqnarray*}
r_{ij}(W(G))&=&\frac{3}{2}+\left[(\frac{3}{5}B^{T}_{1}
+\frac{1}{5}B^{T}_{2})L^{\#}_{G}(\frac{3}{4}B_{1}
+\frac{1}{4}B_{2})\right]_{ii}+\left[(\frac{1}{5}B^{T}_{1}+\frac{3}{5}B^{T}_{2})L^{\#}_{G}
(\frac{1}{4}B_{1}+\frac{3}{4}B_{2})\right]_{jj}\\&&
-[\frac{1}{4}I_{m}+(\frac{3}{5}B^{T}_{1}+\frac{1}{5}B^{T}_{2})
L^{\#}_{G}(\frac{1}{4}B_{1}+\frac{3}{4}B_{2})]_{ij},
\end{eqnarray*}
as stated in $(vi)$.

For any $i\in V_{2}$, $j\in V_{3}$, by Lemma 2.1 and the Equation $(4.2)$, we have
\begin{eqnarray*}
r_{ij}(W(G))&=&\frac{7}{4}+\left[(\frac{2}{5}B^{T}_{1}+\frac{2}{5}B^{T}_{2})L^{\#}_{G}
(\frac{1}{2}B_{1}+\frac{1}{2}B_{2})\right]_{ii}+\left[(\frac{1}{5}B^{T}_{1}+\frac{3}{5}B^{T}_{2})L^{\#}_{G}
(\frac{1}{4}B_{1}+\frac{3}{4}B_{2})\right]_{jj}\\&&
-[\frac{1}{2}I_{m}+(\frac{1}{4}B^{T}_{1}+\frac{1}{4}B^{T}_{2})L^{\#}_{G}
(\frac{1}{4}B_{1}+\frac{3}{4}B_{2})]_{ij},
\end{eqnarray*}
as stated in $(vii)$.

By Lemma 2.4, we have
\begin{eqnarray*}
Kf(W(G))
&=&(n+3m)tr( N)-{\bf{1}^{T}}N{\bf{1}^{T}}\\
&=&(n+3m)\left(\frac{4}{5}tr(L^{\#}_{G})+tr[(\frac{3}{5}B^{T}_{1}
+\frac{1}{5}B^{T}_{2})L^{\#}_{G}(\frac{3}{4}B_{1}+
+\frac{1}{4}B_{2})]+tr[(\frac{2}{5}B^{T}_{1}+\frac{2}{5}B^{T}_{2})
L^{\#}_{G}\right.\\&&\left(\frac{1}{2}B_{1}+\frac{1}{2}B_{2})]+
tr[(\frac{1}{5}B^{T}_{1}+\frac{3}{5}B^{T}_{2})L^{\#}_{G}
(\frac{1}{4}B_{1}+\frac{3}{4}B_{2})]+\frac{5m}{2}\right)-
{\bf{1}^{T}}N{\bf{1}}\\
&=&(n+3m)\left(\frac{4}{5n}Kf(G)
+\frac{61}{100}\left(tr(B^{T}_{1}L^{\#}_{G}B_{1})
+tr(B^{T}_{2}L^{\#}_{G}B_{2})\right)\right.\\&&\left.+\frac{1}{2}
\left(tr(B^{T}_{1}L^{\#}_{G}B_{2})
+tr(B^{T}_{1}L^{\#}_{G}B_{2})\right)+\frac{5m}{2}\right)-
{\bf{1}^{T}}N{\bf{1}}.
\end{eqnarray*}
Since $L^{\#}_{G}1=0$, then
\begin{eqnarray*}
{\bf{1}^{T}}N{\bf{1}^{T}}
&=&1^{T}P_{1}1+1^{T}P_{2}1+1^{T}P_{3}1+
1^{T}Q_{1}1+1^{T}Q_{2}1+1^{T}Q_{3}1+
1^{T}R_{1}1+1^{T}R_{2}1+1^{T}R_{3}1\\&&
2\times1^{T}\frac{3}{4}I_{m}+4\times1^{T}
\frac{1}{2}I_{m}1+2\times1^{T}\frac{1}{4}I_{m}+\times1^{T}I_{m}1\\
&=&
\frac{141}{80}1^{T}B^{T}_{1}L^{\#}_{G}B_{1}1+
\frac{131}{80}1^{T}B^{T}_{1}L^{\#}_{G}B_{2}1
+\frac{133}{80}1^{T}B^{T}_{2}L^{\#}_{G}B_{1}1+
\frac{127}{80}1^{T}B^{T}_{2}L^{\#}_{G}B_{2}1+5m
\end{eqnarray*}

Plugging the above equation into $Kf(W(G))$,
we obtain the required result in $(viii)$.

\section{Conclusion}
In this paper
using the Laplacian generalized inverse
approach, we obtained the resistance distance
and Kirchhoff indices of the quadrilateral graph and the pentagonal graph whenever $G$ is an arbitrary graph.
We can obtain the resistance distance and Kirchhoff
index of the quadrilateral graph and the pentagonal graph
in terms of the resistance distance and kirchhoff index
of the factor graph $G$.

\vskip 0.1in
\noindent{\bf Acknowledgment:}
This work was supported by
the National Natural Science Foundation of China (No.11461020) and
the Research Foundation of the Higher Education
Institutions of Gansu Province, China (2018A-093).


\begin{thebibliography}{99}

\bibitem{KR}D. J. Klein, M. Randi$\acute{c}$, Resistance distance, J. Math. Chem, 12 (1993) 81-95.

\bibitem{BG}R. B. Bapat, S. Gupta, Resistance distance in wheels and fans, Indian J. Pure Appl.Math, 41 (2010) 1-13.


\bibitem{LH}D. Q. Li, Y. P. Hou, The normalized Laplacian spectrum of quadrilateral graphs and its applications, Appl. Math.  Comput. 297 (2017) 180-188.


\bibitem{ZH}H. P. Zhang, Y. J. Yang, C. W. Li, Kirchhoff index of composite graphs, Discrete Appl. Math, 157 (2009) 2918-2927.


\bibitem{BYZhZH}C. J. Bu, B. Yan, X. Q. Zhou, J. Zhou, Resistance distance in subdivision-vertex join and subdivision-edge join of graphs, Linear Algebra Appl., 458 (2014) 454-462.


 \bibitem{LZhB}X. G. Liu, J. Zhou, C. J. Bu, Resistance distance and Kirchhoff index of $R$-vertex join and $R$-edge join of two graphs, Discrete Appl. Math, 187 (2015) 130-139.

  \bibitem{LP}J. B. Liu, X. F. Pan, F. T. Hu, The $\{1\}$-inverse of the Lapalacian
  of subdivision-vertex and suvdivision-edge corona with applications,
  Linear and Multilinear Algebra, 65(2017) 178-191.




 \bibitem{A. Ben-Israel}Ben-Israel, A, Greville, T. N. E., Generalized inverses: theory and applications. 2nd ed.,Springer, New York, 2003.

\bibitem{BSZHW}C. J. Bu, L. Z. Sun,  J. Zhou, Y. M. Wei, A note on block representations of the group inverse of Laplacian matrices, Electron. J. Linear Algebra, 23 (2012) 866-876.



\bibitem{ChZh} H. Y. Cheng, F. J. Zhang, Resistance distance and the normalized Laplacian spectrum, Discrete Appl. Math, 155 (2007) 654-661.








\bibitem{SWZhB}Sun, L. Z., Wang, W. Z., Zhou, J., Bu. C. J., Some results on resistance distances and resistance matrices, Linear and Multilinear Algebra, 63(3) (2015) 523-533.







\bibitem{XG}Xiao, W. J., Gutman, I., Resistance distance and Laplacian spectrum. Theor. Chem. Acc,
110 (2003) 284-289.


\bibitem{YK}Yang, Y. J.,  Klein, D. J., Resistance distance-based graph invariants of subdivisions and triangulations of graphs, Discrete
Appl. Math, 181 (2015) 260-274.



\bibitem{FZZh}Zhang, F. Z, The Schur Complement and Its Applications, Springer-Verlag, New York, 2005.



\end{thebibliography}
\end{document}